# On the Asymptotic Properties of a Certain Class of Goodness-of-Fit Tests Associated with Multinomial Distributions

Sherzod M. Mirakhmedov

V.I. Romanovskiy Institute of Mathematics. Academy of Sciences of Uzbekistan
University str.,46. Tashkent-l00174
e-mail: shmirakhmedov@yahoo.com

**Abstract**. The object of study is the problem of testing for uniformity of the multinomial distribution. We consider tests based on symmetric statistics, defined as the sum of some function of cell-frequencies. Mainly, attention is focused on the class of power divergence statistics, in particular, on the chi-square and log-likelihood ratio statistics. The main issue of the article is to study the asymptotic properties of tests at the concept of an intermediate setting in terms of so called $\alpha$-intermediate asymptotic efficiency due to Ivchenko and Mirakhmedov (1995), when the asymptotic power of tests are bounded away from zero and one, while sequences of alternatives converge to the hypothesis, but not too fast.

## 1. Introduction

The multinomial distribution has received an extensive attention in the literature due to its wide applicability in diverse fields. Assume $n$ particles are allocated into $N$ cells indexed 1 through $N$ at random, successively and independently of each other, the probability of a particle falling into cell with index $l$ is $p_l > 0$, $l = 1,...,N$, $p_1 + ... + p_N = 1$. Let $\eta_l$ be the number of particles in the cell with index $l$ after allocation of all $n$ particles, then $(\eta_1,...,\eta_N)$ has the multinomial distribution $M(n, N, \mathrm{P})$:

$$P\{\eta_1 = m_1,...,\eta_N = m_N\} = \frac{n!}{m_1! \cdot ... \cdot m_N!} p_1^{m_1} \cdot ... \cdot p_N^{m_N},$$

where $\mathrm{P} = (p_1,...,p_N)$ and arbitrary non-negative integer $m_i$s such that $m_1 + ... + m_N = n$.

This probabilistic model arises in statistical mechanics, clinic trails, cryptography, computer theory, and in various fields as an application of the occupancy problem. One can observe that many applications propose some assumption on the cell-probabilities, which need to be tested. In light of this we emphasize that one of the basic tasks in statistics is to ascertain whether a given set of $n$ independent and identically distributed draws come from a given distribution $F_0$, say. This problem is transformed into a problem of fit for a multinomial distribution: support of the given distribution is divided into $N$ mutually exclusive intervals and is counted the number of observations $\eta_m$ arisen in the $m^{\text{th}}$ interval, then the random vector $(\eta_1,...,\eta_N)$ has distribution $M(n,N,\mathrm{P})$. The problem reduces to test the hypothesis $H_0: \mathrm{P}=(p_{01},...,p_{0N})$, where $p_{0m}$ is the probability that a draw has come from $m^{\text{th}}$ interval under the $F_0$. We note that in very important case when $F_0$ is absolutely continuous distribution through probability integral transformation $z \to F_0(z)$ the problem reduces to testing for uniformity $[0,1]$, and then $H_0$ is $\mathrm{P}=(N^{-1},...,N^{-1})$. There seems no need to use statistical methods in the case alternatives far away from the null hypothesis. Therefore, it is of interest to consider the sequences of alternatives approaching $H_0$, viz., $H_{1n}: p_m = p_{0m}(1+\varepsilon_{mn})$, where $p_{01}\varepsilon_{1n}+...+p_{0N}\varepsilon_{Nn}=0$, $\max_{1\le m\le N}|\varepsilon_{mn}| \to 0$ as $n \to \infty$.

The classical tests goodness of fit on the cell-probabilities of cells based on the chi-square type statistics (that is, on the statistics with an asymptotic chi-square distribution), and assume that the number of cells $N$ is fixed. The Pearson's chi-square statistics and the log-likelihood ratio statistics, which are special variants of the power-divergence statistics, introduced by Cressie and Read (1984) are most well-known of them. There is huge literature where interest and results have followed many aspects: the asymptotic distributional and statistical properties and recommendations in applications of power-divergence statistics and its special variants in the case fixed $N$, see Moore (1986), Cressie and Read (1989), and references within. However, the assumption "$N$ is fixed" becomes restrictive in several contexts. Indeed: Mann and Wald (1942) have obtained the relation $N \sim cn^{2/5}$ concerning the optimal choice of the number of groups in chi-square goodness of fit test. Koehler and Larntz (1980) have explored the practical importance of the asymptotic normality results of chi-squared and log-likelihood ratio statistics in case $N$ increases. For the motivation for increasing $N$ associated with "big-data" applications, see Pietrzak et al (2016). At last, a wide class of statistics of interest in the above described multinomial random allocation of particles into cells assumes that the number of cells $N$ increases together with the number of particles $n$, see Kolchin et al (1976), L'ecuyer et al (2002).

In this paper we assume that $N=N(n) \to \infty$ as $n \to \infty$, such that $n/N \to \lambda \in [0,\infty]$ and $n^2/N \to \infty$. We consider statistics of the form $R_N = h_1(\eta_1)+...+h_N(\eta_N)$, where $h_m$s are real-valued functions defined on the non-negative axis. Cohen and Sackrowitz (1975) have proved that the tests based on the statistic $R_N$ are unbiased in testing *uniformity* of multinomial distribution; moreover if $p_{0m}$s are not equal the chi-square test (a special variant of $R_N$) may not be unbiased. Next, if $H_0$ is not uniform or the functions $h_m$ are not the same the tests can't be goodness-of-fit, since their asymptotic power will depend on alternative, and hence they can't distinguish between $H_0$ and all sequences of alternatives of family $H_{1n}$. For the details see Holst (1972) and Ivchenko and Medvedev (1978).

So in this work we consider the testing of $H_0: \mathrm{P} = (N^{-1},...,N^{-1})$ versus sequences of alternatives $H_{1n}: \mathrm{P} = (p_1,...,p_N) \neq (N^{-1},...,N^{-1})$, which *approach* $H_0$ as $n, N \to \infty$ so that

$$\varepsilon(N) = \frac{1}{N}\sum_{m=1}^{N}(Np_m - 1)^2 \to 0, \qquad (1.1)$$

by means of tests based on symmetric statistics of the form

$$S_N^h = \sum_{l=1}^{N} h(\eta_l), \qquad (1.2)$$

where $h$ is a nonlinear real-valued function defined on the non-negative axis. The test based on statistic $S_N^h$ is called $h$-test for short. In what follows we shall assume that the large values of $S_N^h$ reject $H_0$, and $n\lambda_n \to \infty$, as $n \to \infty$, where $\lambda_n = n/N$ is the average of observations per cells. We are interested mainly in the power divergence statistics (PDS) $CR_N(d)$ of Cressie and Read (1984) for which $h(x) = h_d(x)$, where

$$h_d(x) = \frac{2}{d(d+1)} x[(x/\lambda_n)^d - 1], \; d > -1, \; d \neq 0, \text{ else } h_0(x) = 2x\log(x/\lambda_n). \qquad (1.3)$$

One could consider $d \to -1$ and $d < -1$, but this unnecessary in the context of this paper. We emphasize the following important variants of statistics (1.2): the PDSs

$$\chi_N^2 = \lambda_n^{-1}\sum_{m=1}^{N}(\eta_m - \lambda_n)^2, \; \Lambda_N = 2\sum_{m=1}^{N}\eta_m \log(\eta_m/\lambda_n) \text{ and } T_N^2 = 4\sum_{m=1}^{N}(\sqrt{\eta_m} - \sqrt{\lambda_n})^2, \qquad (1.4)$$

which are the chi-square statistic $CR_N(1)$, the log-likelihood ratio statistic $CR_N(0)$ and the Freeman-Tukey statistic $CR_N(-1/2)$, respectively, and the count statistics (CS)

$$\mu_r = \sum_{m=1}^{N}\mathrm{I}\{\eta_m = r\}, \; r \geq 0, \; w_l = \sum_{m=1}^{N}\mathrm{I}\{\eta_m \geq l\}, \; l \geq 1, \text{ and } C_n = \sum_{m=1}^{N}(\eta_m - 1)\mathrm{I}\{\eta_m > 1\}, \qquad (1.5)$$

where $\mathrm{I}\{\cdot\}$ denotes the indicator function, which are respectively, the number of intervals consisting exactly $r$ and at least $r$ observations, and the number of collisions (that is, the number of observations that we observe in intervals already containing observations). These CS have been considered in the literature in various contexts; see, for instance, Kolchin et al (1976), L'ecuyer et al (2002), Khmaladze (2011).

In the classical goodness-of-fit problem of testing for uniformity over [0,1] the most common in the literature is the family of contamination alternatives: $f(x) = 1 + \delta(n)g_n(x)$, where $\delta(n) \to 0$, $\int_0^1 g_n(x)dx = 0$, $0 < \inf_n \|g_n\|_2 \leq \sup_n \|g_n\|_\infty < \infty$, $\|\cdot\|_\infty$ denotes the supremum norm, $\|\cdot\|_2$ is the $L_2[0,1]$ norm. This problem has been intensively studied in the literature; see Inglot (1999), Inglot et al (2019) and references within. We specifically refer to Holst (1972), Gvanceladze and Chibisov (1979), Quine and Robinson (1985), Kallenberg (1985), where the problem was studied in terms of grouped data, and hence the problem is reduced to testing of uniformity of a multinomial distribution against the alternatives

$$H_{1n}: p_m = N^{-1}\left(1+\delta(n)\Delta_{m,n}\right), \ m=1,...,N, \tag{1.6}$$

where $\max_{1\leq m\leq N}|\Delta_{m,n}|\leq C$, $\Delta_{1,n}+...+\Delta_{N,n}=0$ and $N^{-1}(\Delta_{1n}^2+...+\Delta_{Nn}^2)=\Delta^2$ bounded away from zero and infinity. Obviously in this case $\varepsilon(N)=\delta^2(n)\Delta^2$. These alternatives converge to $H_0$ with a rate determined by $\delta(n)$, whereas function $g_n(x)$ (and hence the quantities $\Delta_{mn}$) defines the path along which one goes from the alternative to the hypothesis. For the asymptotic properties of $h$-tests the actual direction of approach to the hypothesis is immaterial (because of its symmetry), but the rate of convergence plays a role. Therefore, it is preferable to present the rate of convergence of alternatives in term of $\varepsilon(N)$, since it does not depend on the direction of convergence.

Our main goal is to study the asymptotic properties of $h$-tests in order to compare them. There are several approaches to the asymptotic comparison of tests which differ by the conditions imposed on the asymptotic behavior of the size, the power and the alternatives, see, for instance, Nikitin (1995). The following three of them are most popular in applications and related to the context of this paper. The Pitman's approach assumes the sequences of alternatives converge to the hypothesis in such a rate that the power for a test of fixed size $\alpha \in (0,1)$, say, has a limit in $(\alpha,1)$. Under Bahadur's approach the power of a test has a limit in $(0,1)$ and alternatives do not approach the hypothesis, then the test is characterized by exponentially decreasing rate of the size. Next, one can consider intermediate between Pitman and Bahadur settings: the power is fixed while the sequences of alternatives approach the hypothesis, but slower than in Pitman setting, the performance of the test is measured by the decreasing rate of the size, but slower than in Bahadur's setting. Note that the rates of convergence of alternatives and the size have to be related.

It should be noted that the concept of an intermediate approach was originally introduced by Kallenberg (1983). Then it was developed for classical problem of testing for uniformity [0,1] by Inglot (1999) and applied to several tests (excluding tests based on grouped data) in series of papers, see Inglot et al (2019) and references within.

Turning to the $h$-tests we refer to Holst (1972), Ivchenko and Medvedev (1978), Quine and Robinson (1985), Mirakhmedov (1987) and Ivchenko and Mirakhmedov (1992) where the Pitman asymptotic efficiency (AE) of $h$-tests was studied completely. In particular, it follows from these works that $h$-tests don't distinguish alternatives approaching $H_0$ at the rate $\varepsilon(N)=o((n\lambda_n)^{-1/2})$. The sequences of alternatives which converge to $H_0$ at the rate $\varepsilon(N)=O((n\lambda_n)^{-1/2})$ form the family of Pitman alternatives, and the chi-square test is optimal within class of $h$-tests in terms of Pitman AE. The asymptotic properties of $h$-tests at the basis of the concept of Bahadur AE (i.e. $\varepsilon(N)$

bounded away from zero) were studied by Ronzhin (1984) who considered the sparse multinomial models, i.e. $\lambda_n \to \lambda \in (0, \infty)$, and a sub-class of *h*-tests satisfying the Cramèr condition:

$$E\exp\{a|h(\xi)|\} < \infty, \exists\, a > 0, \xi \text{ is a Poisson r.v. with parameter } \lambda. \qquad (1.7)$$

Note that the PDS with parameter $d > 0$, in particular the chi-square statistic, does not satisfy the Cramèr condition. We also refer to Quine and Robinson (1985), who showed that for the dense models ($\lambda_n \to \infty$) the chi-square test is inferior to the log-likelihood ratio test in terms of the Bahadur AE, in contrast to the fact that these two tests have the same Pitman AE. The same verdict was confirmed by Kallenberg (1985), who considered "very dense" models when $\lambda_n / N \to \infty$. The family of alternatives (1.1) with $\sqrt{n\lambda_n}\,\varepsilon(N) \to \infty$ forms a family of intermediate alternatives. The asymptotic properties of *h*-tests in the intermediate setting remained less investigated. Only for the sparse multinomial models and sub-class of *h*-tests satisfying the Cramèr condition (1.7) Ivchenko and Mirakhmedov (1995) studied the problem in terms of so-called $\alpha$-intermediate AE ($\alpha$-IAE), which can be interpreted in terms of the decreasing rate of the significance level of *h*-test (see Section 2 below). In particular, their result shows that the chi-square test still optimal (in terms of the $\alpha$-IAE) within class of *h*-tests for the alternatives approaching the hypothesis at the rate $\varepsilon(N) = O(N^{-1/2} \log^{1/2} N)$, but for the alternatives which lie at the "distance" of order $\varepsilon(N) \gg N^{-1/3} \log^{2/3} N$ from the hypothesis the chi-square test is inferior to *h*-tests satisfying the Cramér condition.

Thus the following problems of studying of the intermediate properties of *h*-tests remained open:

- To determine how far from $H_0$ the intermediate alternatives may lie for the chi-square test to retain its asymptotic optimality within the class of *h*-tests;
- To study the intermediate properties of *h*-tests for which the Cramer's condition may not be met;
- To study the intermediate properties of *h*-tests for the "very sparse" and dense models, i.e. when $\lambda_n \to 0$ and $\lambda_n \to \infty$, respectively.

The present work addresses these problems in terms of $\alpha$-IAE. From the results of Section 3, among others, the following significant complements to the result of Ivchenko and Mirakhmedov (1995) follow. Assume the intermediate alternatives are specified by condition $\varepsilon(N) \ll (n\max(1, \lambda_n^2))^{-1/3}$. Then for the sparse models the chi-square test is unique efficient within class of $h_d$-tests and *h*-tests satisfying Cramer condition, whereas the $h_d$-tests and tests based on count statistics $C_n, \mu_r, r = 0, 1, 2$ for the very sparse models, as well as all $h_d$-tests for the dense models have the same efficiency. Next, for the intermediate alternatives which lie at the distance $\varepsilon(N) \gg (n\lambda_n)^{-1/3} \log^{2/3}(N/\lambda_n)$ the chi-square test is much inferior to log-likelihood ratio test if $0 < c \le \lambda_n \ll N$. This fact extends the efficiency properties of these tests in the Bahadur's situation of fixed alternatives presented by Quine and Robinson (1985) to an "adjoin" domain of alternatives approaching $H_0$. Further, the essence of the intermediate setting is that the size tends to 0 as *n* increases, while the asymptotic power, under the underlying sequence of local alternatives, should be non-degenerate. This implies that the rate of convergence of sizes and alternatives have to be linked up. This is shown in fact in the theorems of Section 3.

The rest of the paper is organized as follows. Section 2 provides a concise overview of latest results. The main results of the article are presented in Sections 3 and 4, and the proofs in Section 5. For the reader's convenience, the auxiliary assertions are collected in Appendix.

## 2. Brief survey

In what follows $\lambda_n = n/N$, $\xi \sim Poi(\lambda)$ stands for "a r.v. $\xi$ has Poisson distribution with parameter $\lambda > 0$", $\Phi(u)$ denotes a standard normal distribution function; $c_j$ is a positive constant, may not the same in each its occurrence; *all asymptotic statements are considered as* $n \to \infty$, whenever it is convenient we shall use notation $a_n \ll b_n$ instead of well-known notation $a_n = o(b_n)$.

Asymptotic properties of *h*-tests have been studied by many authors. Let's make a brief review of the latest results. We will distinguish three types of multinomial models: as $n \to \infty$

- the sparse model, when $\lambda_n \to \lambda \in (0,\infty)$,
- the "very sparse" model, when $\lambda_n \to 0$ and $n\lambda_n \to \infty$,
- the dense model, when $\lambda_n \to \infty$.

Let $P_i$, $E_i S_N^h$ and $Var_i S_N^h$ stand for the probability, expectation and variance of $S_N^h$ counted under $H_i, i=0,1$, respectively, $\xi \sim Poi(\lambda_n)$,

$$\tau_n = \lambda_n^{-1} \operatorname{cov}(h(\xi),\xi), \quad g(\xi) = h(\xi) - Eh(\xi) - \tau_n(\xi - \lambda_n),$$

$$\sigma^2(h) = Var\, g(\xi) = Var h(\xi)\left(1 - corr^2\left(h(\xi),\xi\right)\right),$$

$$\rho(h,\lambda_n) = corr\left(h(\xi) - \tau_n \xi, \xi^2 - (2\lambda_n + 1)\xi\right).$$

In what follows we will often refer to the following

**Proposition 2.1.** Assume

$$E|g(\xi)|^3 / \sigma^3(h)\sqrt{N} \to 0, \tag{2.1}$$

and sequences of alternatives $H_{1n}$ satisfy (1.1). Then

$$P_i\left\{S_N^h < u\sigma_i(h)\sqrt{N} + NA_i(h)\right\} = \Phi(u) + o(1), \quad i=0,1, \tag{2.2}$$

and if additionally $\max_m |\varepsilon_{m,n}| = o(\lambda_n^{-1/2})$ then

$$x_N(h) \stackrel{def}{=} \sqrt{N}(A_1(h) - A_0(h))/\sigma_0(h) = \sqrt{n\lambda_n/2}\,\rho(h,\lambda_n)\varepsilon(N)(1+o(1)), \qquad (2.3)$$

where

$$A_i(h) = N^{-1}\sum_{m=1}^{N} E_i h(\xi_m), \quad \xi_m \sim Poi(np_m), \quad \sigma_1^2(h) = \sigma_0^2(h)(1+o(1)) = \sigma^2(h)(1+o(1)). \qquad (2.4)$$

The asymptotical normality result (2.2) follows from Theorem 1 of Mirakhmedov (1992). Proof of (2.3) is given by Mirakhmedov (2022). In particular, from Proposition 2.1 and well-known theorem on convergence of moments (see Theorem 6.14 of Moran (1984)) it follows that $A_i(h)$ and $\sigma_i^2(h)$ are the asymptotic value of $N^{-1}E_i S_N^h$ and $N^{-1}Var_i S_N^h$, respectively.

It follows from (2.2) and (2.3) that there is no power of $h$- tests for the alternatives $H_{1n}$ with $\varepsilon(N) = o((n\lambda_n)^{-1/2})$. Let $\alpha_n(h)$ and $\beta_n(h)$ denote the size and power of $h$-test, respectively. Assume $\alpha_n(h) \to \alpha \in (0,1)$ and $\varepsilon(N) = (n\lambda_n)^{-1/2}$. Then Proposition 2.1 yields

$$\beta_n(h) = \Phi\left(|\rho(h,\lambda_n)|/\sqrt{2} - u_\alpha\right)(1+o(1)), \quad \Phi(-u_\alpha) = \alpha. \qquad (2.5)$$

This implies that the chi-squared test is asymptotically most powerful (AMP) within the family of $h$-tests satisfying (2.1) for all range of the $\lambda_n$, since $|\rho(h,\lambda_n)| \leq 1$ and $|\rho(h,\lambda_n)| = 1$ iff $h(x) = x^2$. The chi-squared test is unique AMP for the sparse multinomial model, but for the very sparse and the dense models the chi-squared test is no longer unique AMP, since there exist test statistics for which $|\rho(h,\lambda_n)| \to 1$ as $\lambda_n \to 0$ or $\lambda_n \to \infty$.

**Remark 2.1**. We emphasize that the condition (2.1) is fulfilled for the sparse models if $E|h(\xi)|^3 < \infty$, and for the very sparse models if $\Delta^2 h(0) \neq 0$, where $\Delta h(x) = h(x+1) - h(x)$. For instance, arbitrary PDS and the CS $C_n, \mu_r, r = 0,1,2$, and $w_l, l = 1,2$ satisfy this condition. But for the dense models the condition (2.1) may impose an additional condition to $\lambda_n$. For instance, (2.1) is fulfilled for PDS (which, remind, include $\chi_N^2$ and $\Lambda_N$ statistics) for arbitrary dense models, while, for example, for CS $\mu_r$, $r \geq 0$ and $C_n$ the (2.1) imposes condition $\lambda_n - \ln N - r \ln \ln N \to -\infty$ and $\lambda_n - \ln N \to -\infty$, respectively. We emphasize also that: if $\lambda_n \to 0$ and statistics such that $\Delta^2 h(0) \neq 0$ then one has, see Lemma 2.2 and 2.3 of Mirakhmedov (2022),

$$\rho(h,\lambda_n) = 1 - \frac{\lambda_n}{6}\left(\frac{\Delta^3 h(0)}{\Delta^2 h(0)}\right)^2 + O(\lambda_n^2), \qquad (2.6)$$

and if $\lambda_n \to \infty$, then for PDS with parameter $d > -1$

$$\rho(h,\lambda_n) = 1 - \frac{(d-1)^2}{6\lambda_n} + O(\lambda_n^{-2}). \tag{2.7}$$

But for the CS (1.5) $\rho(h,\lambda_n) = o(1)$ if $\lambda_n \to \infty$.

**Remark 2.2.** Functional $\rho(h,\lambda_n)$ plays an important role in determining the asymptotic properties of $h$-tests satisfying (2.1). Its sense is clarified by the following fact proved by Mirakhmedov (2022, Lemma 2.1): if $Eh^2(\xi)\xi < \infty$ then

$$\rho(h,\lambda_n) = corr_0(S_N^h, \chi_N^2)(1+o(1)).$$

Further comparison of AMP $h$-tests (for very sparse and dense models) based on the "second order asymptotic efficiency" (SOAE) of $h$-tests w.r.t. the chi-squared test (which is considered as a benchmark). The SOAE based on asymptotic expansion of the power $\beta_n(h)$; it was introduced and studied by Ivchenko and Mirakhmedov (1992). Set $t_\alpha = 2^{-1/2} - u_\alpha$, see (2.5).

**Definition 2.1.** An AMP $h$-test is called SOAE wrt the chi-square test, when $\lambda_n \to 0$ or $\lambda_n \to \infty$, if in asymptotic expansion of its power, viz. $\beta_n(h) = \Phi(t_\alpha) + \vartheta_N^h(1+o(1))$, the term $\vartheta_N^h$ asymptotically coincides with similar term in asymptotic expansion of the power of chi-square test.

The following statements have been proved by Ivchenko and Mirakhmedov(1992) and Mirakhmedov et al (2014, p.738).

*For very sparse model*: the $h$-tests such that $\Delta^2 h(0) \neq 0$ are constitute a sub-family of AMP $h$-tests; for the power of $h$-test the following type asymptotic expansion hold

$$\beta_n(h) = \Phi(t_\alpha) - \phi_1^h(t_\alpha)\lambda_n(1+o(1)) + \phi_2^h(t_\alpha)(n\lambda_n)^{-1/2}(1+o(1)),$$

at that for the chi-square test

$$\phi_1^h(t_\alpha) = 0, \quad \phi_2^h(t_\alpha) = \frac{1}{\sqrt{2\pi}} e^{-t_\alpha^2/2}\left(\frac{1-t_\alpha^2}{3\sqrt{2}} + \frac{t_\alpha}{2} + \sqrt{2}S_1\left(t_\alpha\sqrt{\frac{n\lambda_n}{2}} + \frac{n}{2}\right)\right),$$

where $S_1(u) = u - \lfloor u \rfloor + 1/2$, $\lfloor u \rfloor$ is the integer part of $u$. From these it follows that if $N = o(n^{4/3})$, i.e. $n^{-1/3} = o(\lambda_n)$, then no SOAE test exists within family of $h$-tests; but if $n = O(N^{3/4})$, then there exist $h$-tests with $\Delta^3 h(0) \neq 0$, which are SOAE, the CS (1.6) are examples of such tests. Thus, the SOAE $h$-tests exist in multinomial models as sparse as $\lambda_n = O(N^{-1/4})$.

*The dense model* turned out to be more complicated. It is only shown that the PDS constitute the family of AMP $h$-tests; whereas the CSs can't generate AMP test. Next, the chi-squared test still optimal within family of PDS in sense of SOAE if $n = o(N^{3/2})$. On the basis of formal asymptotic expansion of power functions it has been concluded that the log-likelihood test can't be SOAE for $n \geq N^{3/2}$.

Let's turn to an "intermediate properties" of $h$-tests. First we make some comments to clarify the meaning of "intermediate properties" of $h$-tests, which rises in context of comparison of two tests.

As we pointed out in introduction, the concepts for comparing the performance of two sequences of statistical tests for a given hypothesis testing problem are differ by the conditions imposed on the asymptotic behavior $\alpha_n(h)$, $\beta_n(h)$ and alternatives $H_{1n}$ (by a condition on the $\varepsilon(N)$ in our case). The conditions imposed on two of them provide a condition for the third. The Pitman's approach assumes that $\alpha_n(h) \to \alpha > 0$ and a sequence of alternatives converge to the hypothesis at the rate necessary to $\beta_n(h) \to \beta \in (\alpha, 1)$. In our problem this Pitman sub-family of alternatives is $H_{1n}$ where $\varepsilon(N) = (n\lambda_n)^{-1/2}$. It follows from Proposition 2.1 and (2.3) that the Pitman asymptotical efficiency of $h$-test wrt $g$-test, see Nikitin (1995, p. 26), for Pitman sub-family of alternatives is $PE(S_N^h, S_N^g) = \lim \rho^2(h, \lambda_n) / \rho^2(g, \lambda_n)$. In particular, $PE(\chi_N^2, S_N^h) = \lim \rho^{-2}(h, \lambda_n) \geq 1$, that is within the class of $h$-tests the chi-square test is the most efficient in Pitman's sense.

Another "extreme" setting assumes that $\beta_n(h) \to \beta \in (0,1)$ and the alternatives $H_1$ do not approach the $H_0$, for our problem this means $\varepsilon(N)$ is a constant. This approach is at the basis of the concept of Bahadur AE; it was developed by Ronzhin (1984), who considered the sparse multinomial models with "small samples", when $\max_m np_m \leq c$, and a sub-class of $h$-tests satisfying the Cramèr condition (1.7). This condition means that $h(x) = O(x \log x)$, and hence significantly limits the class of $h$-tests. Set $\xi(z) \sim Poi(z)$, $\psi_h(t,z) = E \exp\{th(\xi(z))\}$, and $c_h(t) = \lambda \log \lambda - \lambda + z_h(t) - \lambda \log z_h(t) + \log \psi_h(t, z_h(t))$, where $z_h(t)$ is a solution of the equation $E(\xi(z) - \lambda) \exp\{th(\xi(z))\} = 0$. Remain $h$-test rejects $H_0$ for large values of $S_N^h$. Assume that $\lim_{N \to \infty}(A_1(h) - A_0(h)) > 0$, see (2.4), this means that the alternatives do not approach $H_0$. Due to Ronzhin (1984) the Bahadur's exact slope of $h$-test is

$$-\lim_{N \to \infty} N^{-1} \log P_0\{S_N^h > NA_1(h)\} = t_0 \tilde{A}_1(h) - c_h(t_0) := J(h, H_1),$$

where $\tilde{A}_1(h)$ is limit of $A_1(h)$ as $N \to \infty$, $t_0$ is a solution of equation $c_h'(t_0) = \tilde{A}_1(h)$, $c_h'(t)$ is a derivative of $c_h(t)$ w.r.t. $t$, and this limit specifies the Bahadur efficiency of $h$-test. Let $\Upsilon_m$ denotes a sub-class of symmetric statistics generated by function $h(x, a_0, ..., a_m) = a_0 I\{x=0\} + ... + a_m I\{x=m\}$, where an integer $m \geq 0$. Ronzhin showed that if $N^{-1} E_1 \mu_r \to b_r$, $r = 0, ..., m$, and $(b_0, ..., b_m) \neq (\pi_0(\lambda), ..., \pi_m(\lambda))$, where $\mu_r$ is CS defined in (1.5), $\pi_r(\lambda) = \lambda^r e^{-\lambda} / r!$, then the Bahadur efficient $h$-test within class $\Upsilon_m$ is generated by function $\tilde{h}(x, a_0, ..., a_m)$, where

$$a_r = \log\frac{b_r}{\pi_r(\omega)} - \log\frac{1-b_0-...-b_m}{1-\pi_0(\omega)-...-\pi_m(\omega)}, \quad r=0,...,m,$$

and $\omega$ is the solution to the equation

$$\omega\frac{1-b_0-...-b_m}{\lambda-1b_1-...-mb_m} = \frac{1-\pi_0(\omega)-...-\pi_m(\omega)}{1-\pi_0(\omega)-...-\pi_{m-1}(\omega)}.$$

Note that the optimal function $\tilde{h}$ depend in this case on the alternatives, through vector $(b_0,...,b_m)$, and for this function

$$J(\tilde{h}, H_1) = \lambda - \omega + \lambda\log\frac{\omega}{\alpha} + I_m(b,\pi),$$

where $I_m(b,\pi)$ is the Kullback –Sanov's distance

$$I_m(b,\pi) = \sum_{r=0}^{m} b_r\log\frac{b_r}{\pi_r(\omega)} + (1-b_0-...-b_m)\log\frac{1-b_0-...-b_m}{1-\pi_0(\omega)-...-\pi_{m-1}(\omega)}.$$

In particularly, Bahadur efficiency of the statistic $\mu_0$ (the case $m=0$) is equal to

$$J(b_0) = \lambda - \omega + \lambda\log\frac{\omega}{\alpha} + I_0(b_0, e^{-\lambda}).$$

An intermediate approach, which lies somewhere between Pitman's and Bahadur's settings, should assume that the asymptotic power is bounded away from zero and one, while the alternatives converge to the hypothesis, but slower than in the Pitman's approach. According to Ivchenko and Mirakhmedov (1995) the performance of $h$-tests in this situation will be measured similarly to Bahadur's setting by the asymptotic value of

$$e_N^\alpha(S_N^h) = -\log P_0\left\{S_N^h > NA_1(h)\right\} = -\log P_0\left\{\tilde{S}_N^h > x_N(h)\right\}, \tag{2.8}$$

where $\tilde{S}_N^h = (S_N^h - NA_0(h))/\sqrt{N}\sigma_0(h)$, $x_N(h) = \sqrt{N}(A_1(h) - A_0(h))/\sigma_0(h)$, as above. Note that $x_N(h)$ coincides with "efficacy" of statistic $S_N^h$, the notion introduced by Freser (1957) for finding of Pitman asymptotic relative efficiency of two tests.

It is seen that the intermediate setting must impose two conditions, in notation (1.1) and (2.3), $\sqrt{n\lambda_n}\,\varepsilon(N)\to\infty$ (in contrast to Pitman's case) and $A_1(h) - A_0(h) \to 0$ (in contrast to Bahadur's situation), which, due to (2.8) and (2.3), specifies the family of intermediate alternatives $H_{1n}$ such that

$$\sqrt{n\lambda_n}\,\varepsilon(N)\to\infty \text{ and } \varepsilon(N) = o(\lambda_n^{-1}). \tag{2.9}$$

The described situation gives rise to the concept of intermediate asymptotic efficiency, which is called $\alpha$-IAE for short. For any fixed $z\in R$ define the significance level $\alpha_n(z,h)$ of the $h$-test corresponding to the critical region $\{\tilde{S}_N^h > x_N(h) + z\}$. Due to Proposition 2.1 statistic $\tilde{S}_N^h$ is bounded in probability under $H_0$, and hence $\alpha_n(z,h) \to 0$ for every sequences of intermediate alternatives (2.9) since $x_N(h) \to \infty$, but this convergence is not exponentially fast, as it is in Bahadur setting. Next, by Proposition 2.1 one can show that under alternatives (2.9) the asymptotical power of this $h$-tests for arbitrary $z$ is equal to $\Phi(-z)$, hence bounded away from 0 and 1. On the other hand the quantity $e_N^\alpha(S_N^h)$ asymptotically coincides with $-\log \alpha_n(z,h)$. Hence, $h$-test that has largest asymptotical value of $e_N^\alpha(h)$ should be considered as efficient within class of $h$-tests.

For the alternatives (2.9) $x_N(h) \to \infty$, $x_N(h) = o(\sqrt{N})$, so asymptotic analysis of $e_N^\alpha(h)$ based on the probabilities of large deviation results under the hypothesis, the order of large deviation depend on the rate of $\sqrt{n\lambda_n}\varepsilon(N) \to \infty$, since (2.3). The "distance" $\varepsilon(N)$ determines the pertaining range of large deviation for the statistic $S_N^h$; this effect is common for all $h$-tests, and hence the $\alpha-\mathrm{IAE}$ of various $h$-tests differ through the functional $|\rho(h,\lambda_n)|$, similarly to Pitman AE.

*The statement of Ivchenko and Mirakhmedov (1995) in our notation reads as follows.*

**Theorem 2.1**. Let $\lambda_n \to \lambda \in (0,\infty)$, the alternative $H_{1n}$ specified by (2.9). Then

$$\frac{e_N^\alpha(S_N^h)}{n\lambda_n \varepsilon^2(N)} = \frac{1}{4}\rho^2(h,\lambda_n)(1+o(1)), \qquad (2.10)$$

provided either

(i) $E|h(\xi)|^{2+\delta} < \infty$, some $\delta > 0$, and $\varepsilon(N) = O(N^{-1/2}\log^{1/2} N)$, or

(ii) $E\exp\{a|h(\xi)|\} < \infty$, some $a > 0$, $\xi \sim Poi(\lambda)$.

This implies that for the sparse models the $\alpha-\mathrm{IAE}$ of $h$-test is determined by the functional $\rho(h,\lambda_n)$, and hence a Pitman efficient $h$-test still optimal in terms of $\alpha-\mathrm{IAE}$, as long as the conditions (i) and (ii) are fulfilled. In particular, the chi-squared statistic satisfy condition (i) but not (ii), so chi-squared test is optimal within class of $h$-tests in terms of $\alpha-\mathrm{IAE}$ if alternatives (2.9) such that $\varepsilon(N) = O(N^{-1/2}\log^{1/2} N)$, whereas for the alternatives at distant $\varepsilon(N) \gg N^{-1/2}\log^{1/2} N$ the chi-squared test is much inferior to $h$-tests satisfying the Cramèr condition, in particularly to tests based on $\Lambda_N$, $T_N^2$, see (1.4), and CS (1.5). Next, the optimality of $h$-tests for the alternatives (2.9) with $\varepsilon(N) \gg N^{-1/2}\log^{1/2} N$ can be deduced only for some sub-class of $h$-tests satisfying the Cramèr condition. For instance, for the sub-class $\{h: h(x)=0, \text{ for } x>m\}$, fixed $m \geq 1$, the condition (ii) is fulfilled, and in terms of $\alpha-\mathrm{IAE}$ the optimal $h$-test within this sub-class coincides with the optimal test in Pitman AE sense presented by Victorova and Chistyakov (1966).

**3. New results**

From statistical point of view, in particular the problem of hypothesis testing, the most important flexible subfamily of (1.2) is the family of PDS $CR_N(d)$, which generates $h_d$-tests, see (1.3). Now on symbol $\Im_{alt}$ stands for the family of alternatives (2.9); $\varepsilon(N) \ll a_n$ stands for the "subfamily of $\Im_{alt}$ specified by condition $\varepsilon(N) \ll a_n$".

**Theorem 3.1**. Let $\lambda_n \to \lambda > 0$. If either

(i) $-1 < d \leq 0$ and the family of alternatives is $\Im_{alt}$, or

(ii) integer $d \geq 1$ and $\varepsilon(N) \ll n^{-d/(1+2d)}$, or

(iii) non-integer $d > 0$ and $\varepsilon(N) \ll \min(n^{-3/8}, n^{-d^*/(1+2d^*)})$, where $d^* = \max(1, d)$. Then

$$\frac{e_N^\alpha(S_N^{h_d})}{n\lambda_n \varepsilon^2(N)} = \frac{1}{4}\rho^2(h_d, \lambda_n)(1 + o(1)). \tag{3.1}$$

**Remarks 3.1**. Part (i), in fact, consists a PDS satisfying Cramèr condition (1.7), for example log-likelihood statistic and the Freeman-Tukey statistic are included in part (i). The chi-square statistic is included in part (ii), where $\varepsilon(N) \ll n^{-1/3}$. The PDS $CR_N(3/2)$ and $CR_N(2/3)$ recommended by Cressie and Read (1984, p.463) are examples of part (iii) where for these cases $\varepsilon(N) \ll N^{-3/8}$.

**Theorem 3.2.** Let $\lambda_n \to 0$, $n\lambda_n^3 \to \infty$. If either

(i) $d \in (-1, 0]$ or integer $d \geq 1$, while $\varepsilon(N) \ll (nd^*\lambda_n^{d^*-1})^{-1/(2d^*+1)}$, where $d^* = \max(1, d)$, or

(ii) non-integer $d > 0$ and $\varepsilon(N) \ll \min\left((n\lambda_n^{4/3})^{-3/8}, (n^{d^*}\lambda_n^{d^*-1})^{-1/(2d^*+1)}\right)$. Then equation (3.1) is valid.

Note that, see Mirakhmedov (2022, Lemma 2.2, if $\lambda_n \to 0$ then

$$\rho(h_d, \lambda_n) = 1 - \frac{3(3^d - 2^{d+1} + 1)^2}{8(2^d - 1)^2}\lambda_n + O(\lambda_n^2), d \neq 0,$$

$$\rho(h_0, \lambda_n) = 1 - \frac{3}{8}\left(\frac{\ln 3/4}{\ln 2}\right)^2 \lambda_n + O(\lambda_n^2).$$

Therefore, it follows from Theorem 3.2 that all $h_d$-tests have the same $\alpha - \text{IAE}$. The parameter $d$ of PDS affects to the family of alternatives for which $h_d$-test is applied; let's denote this family by symbol $\Im_{alt}(d)$. Then $\Im_{alt}(d) \subseteq \Im_{alt}(1)$ for all $d \in (-1, 0]$, and $d > 3/2$ if $\lambda_n \geq n^{-1/3}$, i.e. $N^{3/4} \leq n \ll N$. Hence, w.r.t. these $h_d$-tests the chi-square test ($h_1$-test) is preferable. The same verdict is hold for non-integer $d \in (0, 3/2]$ if $\lambda_n \geq n^{-1/4}$, i.e. $N^{4/5} \ll n \ll N$. Notice that $\Im_{alt}(1)$ is the sub-family of intermediate alternatives such that $\varepsilon(N) \ll n^{-1/3}$.

**Theorem 3.3**. Let $\lambda_n \to \infty$, then equation (3.1) is valid for every $d > -1$ and $\varepsilon(N) \ll (n\lambda_n^2)^{-1/3}$.

By Lemma 2.3 of Miakhmedov (2022) we have $\rho(h_d, \lambda_n) = 1 - (d-1)^2 \lambda_n^{-1} + O(\lambda_n^{-2})$, if $\lambda_n \to \infty$. Thus, Theorem 3.3 says that for dense models all $h_d$-tests with $d > -1$ have the same $\alpha - \text{IAE}$ within sub-family of $\Im_{alt}$ such that $\varepsilon(N) \ll (n\lambda_n^2)^{-1/3}$.

The intermediate properties of chi-square and log-likelihood tests are collected in the following theorems.

**Theorem 3.4.** If either

(i) $n\lambda_n^3 \to \infty$ and $\varepsilon(N) \ll \left(n\max(1, \lambda_n^2)\right)^{-1/3}$; or

(ii) for each $\gamma \in (1/4, 1/3]$ the parameter $N \in N_\gamma = \{N : n^{(1-3\gamma)/(1-2\gamma)} \ll N \ll n^{3(1-2\gamma)/4(1-\gamma)}\}$ and $\varepsilon(N) = (n\lambda_n^2)^{-\gamma}$, then

$$\frac{e_N^\alpha(\chi_N^2)}{n\lambda_n \varepsilon^2(N)} = \frac{1}{4}(1 + o(1));$$

(iii) If

$$\varepsilon(N) = (n\lambda_n)^{-1/3} \omega_n^{2/3}, \text{ where } \max(1, \log(N^2/n)) \ll \omega_n \ll \sqrt{n\lambda_n}, \quad (3.2)$$

then

$$\frac{e_N^\alpha(\chi_N^2)}{n\lambda_n \varepsilon^2(N)} = o(1);$$

**Remarks 3.3.** The condition $N \in N_\gamma$ implies $N \ll \sqrt{n}$ for all $\gamma \in (1/4, 1/3]$. Therefore, parts (i) and (ii) together cover the case $\varepsilon(N) \ll (n\lambda_n^2)^{-1/4}$ if $N \ll \sqrt{n}$. Let $\Im'_{alt}$ and $\Im''_{alt}$ stand for the sub-families of $\Im_{alt}$ such that $\varepsilon(N) \geq (n\lambda_n)^{-1/3}$ and $\varepsilon(N) \geq (n\lambda_n^2)^{-1/4}$, respectively. The family of alternatives satisfying condition (3.2) coincides with $\Im'_{alt}$ if $N = c\sqrt{n}$. Obviously $\Im'_{alt} = \Im''_{alt}$ when $N = \sqrt{n}$; if $N \ll \sqrt{n}$ then $\Im'_{alt} \subset \Im''_{alt}$, and $\Im'_{alt} \supset \Im''_{alt}$ if $N \gg \sqrt{n}$. Hence part (iii) consists wider family of alternatives when $N \gg \sqrt{n}$; part (i) also includes the case $N \gg \sqrt{n}$, when $\varepsilon(N) \ll (n\lambda_n^2)^{-1/3}$. Condition (3.2) can be replaced by stronger condition $\varepsilon(N) \gg (n\lambda_n)^{-1/3} \log^{2/3} N$, which coincides with (3.2) for a sparse models. Thus *remain open* the $\alpha - \text{IAE}$ of chi-square test for the sub-family of alternatives such that $(n\lambda_n^2)^{-1/3} \leq \varepsilon(N) = O\left((n\lambda_n)^{-1/3} \log^{2/3}(N^2/n)\right)$ when $\sqrt{n} \leq N$, and for $\varepsilon(N) \geq (n\lambda_n^2)^{-1/4}$ in the situation $N \ll \sqrt{n}$.

**Theorem 3.5.** If either

(i) $\lambda_n \to \lambda > 0$ and the family of alternatives is $\Im_{alt}$, or

(ii) $n\lambda_n^3 \to \infty$ and $\varepsilon(N) \ll (n\max(1, \lambda_n^2))^{-1/3}$, or

(iii) $\sqrt{n} \ll N \ll n$, for each $\gamma \in (0, 1/3]$ parameter $N \in N_\gamma^{**} = \{n^{(1-2\gamma)/(1-\gamma)} \ll N \ll n^{(5-8\gamma)/(5-4\gamma)}\}$ and $\varepsilon(N) = (n\lambda_n)^{-\gamma}$, or

(iv) $N \ll \sqrt{n}$, for each $\gamma \in (1/4, 1/3]$ parameter $N \in N_\gamma^* = \{N : n^{(1-3\gamma)/(1-2\gamma)} \ll N \ll n^{(5-12\gamma)/(5-6\gamma)}\}$ and $\varepsilon(N) = (n\lambda_n^2)^{-\gamma}$,

then

$$\frac{e_N^\alpha(\Lambda_N)}{n\lambda_n \varepsilon^2(N)} = \frac{1}{4}(1 + o(1)).$$

**Remarks 3.4.** Parts (ii) and (iv) of Theorem 3.5 together cover the $\alpha - AIE$ properties of log-likelihood test for $\varepsilon(N) \ll (n\lambda_n^2)^{-1/4}$, when $N \ll \sqrt{n}$. Note that sets $N_\gamma^{**}$, $\gamma \in (0, 1/3]$, all together cover the set $\{N : \sqrt{n} \ll N \ll n\}$ completely. The case $N = cn$ is included in part (i). Next, the sub-families of alternatives of part (iii) all together form the sub-family of $\Im_{alt}$, which is specified by condition $\varepsilon(N) \geq (n\lambda_n)^{-1/3}$. Part (ii) also includes the case $\sqrt{n} \ll N$, but $\varepsilon(N) \ll (n\lambda_n^2)^{-1/3}$. So it is *remain open* the $\alpha - AIE$ of log-likelihood test for $(n\lambda_n^2)^{-1/3} \leq \varepsilon(N) \ll (n\lambda_n)^{-1/3}$ in the case $\sqrt{n} \ll N$, and for $\varepsilon(N) \geq (n\lambda_n^2)^{-1/4}$ in the situation $N \ll \sqrt{n}$.

Let's consider now the CS (1.5). First we make the following comments. For CS the Cramer condition is fulfilled, hence part (ii) of Theorem 2.1 is applicable for the sparse models. For CS (1.5) defined on dense models we have $\rho(h, \lambda_n) \to 0$, hence their $\alpha - IAE$ in this case is inferior wrt any PDS. At last, for the CS $\mu_r$ defined on the very sparse models again $\rho(h, \lambda_n) \to 0$ if $r \geq 3$, so these CS are not of interest also. Thus it is reasonable to consider CS $\mu_r$, $r = 0, 1, 2$ only, since $\rho(h, \lambda_n) \to 1$ by (2.6).

**Theorem 3.6.** If $\lambda_n \to 0$ and $n^{1/6}\lambda_n \to \infty$, then

$$e_n^\alpha(\mu_r) = 4^{-1} n\lambda_n \varepsilon^2(N)(1 + o(1)), \quad r = 0, 1, 2$$

for every sequence of alternatives of $\Im_{alt}$ such that $\varepsilon(N) \ll N^{-1/3}$. Alike assertion is true for the CS $C_n$.

**4. Concluding remarks**

Theorems 3.1-3.6 allow us to make the following concluding remarks. Recall notation $\tilde{S}_N^h = (S_N^h - NA_0(h))/\sigma_0(h)\sqrt{N}$, and that $x_N(h) = \sqrt{N}(A_1(h) - A_0(h))/\sigma_0(h) \sim \sqrt{n\lambda_n/2}\varepsilon(N)\rho(h, \lambda_n)$ for the family $\Im_{alt}$. For any fixed $z \in R$ define the significance level $\alpha_n(z, h)$ of the $h$-test corresponding to the critical region $\{\tilde{S}_N^h > x_N(h) + z\}$, viz., $\alpha_n(z, h) = P_0\{\tilde{S}_N^h > z + x_N(h)\}$. Note that $-\log \alpha_n(z,h)$ asymptotically coincides with $e_N^\alpha(S_N^h)$, for any $z \in \mathbb{R}$, since $x_N(h) \to \infty$ as $n \to \infty$. Hence, in fact, Theorems 3.1-3.6 establish a relationship between the significance level $\alpha_n(z,h_d)$ and the rate of convergence of alternatives to $H_0$; for instance from equation (3.1) we obtain $\alpha_n(z,h_d) \sim \exp\{-n\lambda_n \varepsilon^2(N) \rho^2(h_d, \lambda_n)/4\}$. In particularly, for the chi-square test and intermediate alternatives with $\varepsilon(N) = (n\lambda_n^2)^{-\gamma}$ and $N \in \mathbb{N}_\gamma = \{N : n^{(1-3\gamma)/(1-2\gamma)} \ll N \ll n^{3(1-2\gamma)/4(1-\gamma)}\}$, $\gamma \in (1/4, 1/3]$ from Theorem 3.5 (ii) we have

$$\exp\{-c\sqrt{n}\} \ll \alpha_n(z,h_1) = \exp\{-4^{-1} n^{1-2\gamma} \lambda_n^{1-4\gamma}(1+o(1))\} \le \exp\{-N^{1/3}/4\}.$$

Note yet, the exact rate of convergence of $\alpha_n(z, h_1)$ depend on the relation between the sample size $n$ and the number of groups (cells) $N$ for each $\gamma$.

Further, by Proposition 2.1 it is seen that for every sequence of alternatives of $\Im_{alt}$

(i) The asymptotical power of arbitrary $h$-test is bounded away from 0 and 1,

(ii) For arbitrary small $\varepsilon > 0$ it holds $P_1\{|x_N^{-1}(h)\tilde{S}_N^h - 1| \ge \varepsilon\} \to 0$ as $n \to \infty$.

Therefore, $e_N^\alpha(S_N^h)$ asymptotically coincides with intermediate slope of $h$-test, in accord with Inglot (1999, p.491). In view of these reasons we introduce the following definition of asymptotic relative efficiency of $h$-test with respect to $g$-test, which is analogue of asymptotic intermediate efficiency in weak sense of Inglot (1999).

**Definition 3.1**. If for every sequence of alternatives of $\Im \subseteq \Im_{alt}$ it holds

$$\lim_{n \to \infty} \frac{e_N^\alpha(S_N^h)}{e_N^\alpha(S_N^g)} = e \ge 0, \quad (3.3)$$

we say that $e := e(h,g)$ is $\alpha$-IAE of $h$-test with respect to $g$-test for the family $\Im$.

We will take the chi-square test ($h_1$-test) as a benchmark procedure and compare other $h$-tests to it. Let $\{h_d\}$ and $\{C\}$ stands for the class of $h_d$-tests, $d > -1$, and of $h$-tests satisfying the Cramèr condition, respectively. Recall that $h_{-1/2}$-test and $h_0$-test are tests based on $T_N^2$ and $\Lambda_N$ statistics, respectively.

**The very sparse models**: Let $\lambda_n \to 0$. Applying Theorem 3.2, Theorem 3.4 (i), Theorem 3.5(ii) and Theorem 3.6 we conclude that $e(h_1, h) = 1$, where $h(x) \in \{h_d\}$ and $\varepsilon(N) \ll (n^{d^*} \lambda_n^{d^*-1})^{-1/(2d^*+1)}$, $d^* = \max(1,d)$, if $n\lambda_n^3 \to \infty$, and also $h(x) = I\{x = r\}$, $r = 0,1,2$, or $h(x) = (x-1)I\{x > 1\}$ if additionally $n\lambda_n^6 \to \infty$ and $\varepsilon(N) \ll N^{-1/3}$. That is the tests based on PDS with any parameter $d > -1$, and CS $\mu_r$, $r = 0,1,2$, and $C_n$ are equally efficient in sense of $\alpha$-IAE for indicated family of intermediate alternatives. In particular, $e(h_1, h_0) = 1$ for $\varepsilon(N) \ll n^{-1/3}$ if $n\lambda_n^3 \to \infty$.

**The sparse models**: $\lambda_n \to \lambda \in (0,\infty)$. From Theorem 3.1 it follows that

$$e(h,g) = \rho^2(h,\lambda)/\rho^2(g,\lambda)(1+o(1)) \qquad (3.4)$$

(i) for $\varepsilon(N) \ll n^{-1/3}$ and $h, g \in \{h_d\} \cup \{C\}$, and (ii) for the $\mathfrak{I}_{alt}$ completely and $h, g \in \{C\}$. These facts extend the Pitman AE properties of these tests presented by Holst (1972) and Ivchenko and Medvedev (1978) to an "adjoin" domain of intermediate alternatives specified by condition $\varepsilon(N) \ll n^{-1/3}$ and $\mathfrak{I}_{alt}$ for indicated class of tests, respectively. Further, for the $\varepsilon(N) \ll n^{-1/3}$ we have $e(h_1, h) = \rho^{-2}(h,\lambda) > 1$, but $e(h_1, h) = 0$ for the $\varepsilon(N) \gg n^{-1/3} \log^{2/3} n$. That is for $\varepsilon(N) \ll n^{-1/3}$ the chi-square test is unique optimal within the class $\{h_d\} \cup \{C\}$, whereas for $\varepsilon(N) \gg n^{-1/3} \log^{2/3} n$ is much inferior wrt tests satisfying the Cramèr condition (1.7) (for instance to tests $h_{-1/2}, h_0$, and tests based on CS (1.5)). This is significant extension of result of Ivchenko and Mirakhmedov (1995).

**The dense model** $\lambda_n \to \infty$: First of all we emphasize that $e(h_{d_1}, h_{d_2}) = 1$ for any $d_1, d_2$ and $\varepsilon(N) \ll (n\lambda_n^2)^{-1/3}$ without any restriction to the increasing rate of $\lambda_n$. It is interesting that for the alternatives at a distant $\varepsilon(N) \geq (n\lambda_n^2)^{-1/3}$ from the hypothesis the $\alpha$-AIE of $h_d$-tests are differ for $\sqrt{n} \ll N$ and $\sqrt{n} \gg N$. Indeed. Let $\sqrt{n} \ll N$, i.e. $\lambda_n \ll \sqrt{n}$. Then it follows from Theorem 3.4(iii) and Theorem 3.5(iii) that $e(h_1, h_0) = 0$ for $\varepsilon(N) \gg (n\lambda_n)^{-1/3} \log^{2/3}(N/\lambda_n)$. This fact extends the efficiency properties of chi-square and log-likelihood ratio tests in the Bahadur's situation of fixed alternatives to an "adjoining" domain of alternatives approaching $H_0$. However, if $N \ll \sqrt{n}$ then from Theorem 3.4 (ii) and Theorem 3.5 (iv) we obtain that $e(h_1, h_0) = 1$ for alternatives at distant $\varepsilon(N) = (n\lambda_n^2)^{-\gamma}$ and $N \in \mathrm{N}_\gamma \cap \mathrm{N}_\gamma^*, \gamma \in (1/4, 1/3]$ (in fact $(n\lambda_n^2)^{-1/3} \leq \varepsilon(N) \ll (n\lambda_n^2)^{-1/4}$). These facts extend the Pitman efficiency properties of the chi-square and log-likelihood ratio tests to an "adjoining" domain of intermediate alternatives specified by the condition $\varepsilon(N) \ll (n\lambda_n^2)^{-1/4}$, when $N \ll \sqrt{n}$.

Further, actually the equality (3.4) holds for very sparse and dense models also (for a suitable family of alternatives from the corresponding theorems) with notice that for these models $\rho(h_d, \lambda_n) \to 1$. So it is seen that $\alpha$-AIE of $h$-test depend on the asymptotic behavior of the parameter $\lambda_n$ and $|\rho(h,\lambda)|$, the correlation coefficient between the test statistic $S_N^h$ and the chi-square statistic; so a statistic that is more correlated with the chi-square statistic should be considered preferable. For the PDS in Table 1 the values of $|\rho(h_d, \lambda)|$ are presented for various $\lambda$ and $d > -1$.

Table 1

| d | $\lambda$ | | | | | | | | | |
|---|---|---|---|---|---|---|---|---|---|---|
|  | 0.05 | 0.1 | 0.5 | 1.0 | 1.5 | 2.0 | 3.0 | 10 | 20 | 50 |
| -2/3 | 0.9933 | 0.9838 | 0.9400 | 0.8768 | 0.8314 | 0.7811 | 0.7266 | 0.9257 | 0.9740 | 0.9900 |
| -1/2 | 0.9942 | 0.9838 | 0.9402 | 0.8909 | 0.8545 | 0.8321 | 0.8001 | 0.9480 | 0.9803 | 0.9920 |
| -1/3 | 0.9950 | 0.9839 | 0.9620 | 0.9192 | 0.89891 | 0.8743 | 0.8573 | 0.9615 | 0.9834 | 0.9940 |
| 0 | 0.9970 | 0.9940 | 0.9720 | 0.9525 | 0.9400 | 0.9350 | 0.9369 | 0.9793 | 0.9897 | 0.9960 |
| 1/3 | 0.9983 | 0.9840 | 0.9845 | 0.9758 | 0.9699 | 0.9714 | 0.9797 | 0.9928 | 0.9961 | 0.9980 |
| 1/2 | 0.9989 | 0.9979 | 0.9900 | 0.9898 | 0.9815 | 0.9791 | 0.9879 | 0.9972 | 0.9993 | 0.9985 |
| 2/3 | 0.9999 | 0.9924 | 0.9901 | 0.9900 | 0.9930 | 0.9945 | 0.9961 | 0.9977 | 0.9996 | 0.9990 |
| 1 | 1.00 | 1.00 | 1.00 | 1.00 | 1.00 | 1.00 | 1.00 | 1.00 | 1.00 | 1.00 |
| 3/2 | 0.9984 | 0.9844 | 0.9900 | 0.9901 | 0.9930 | 0.9925 | 0.9879 | 0.9977 | 0.9997 | 0.9989 |
| 2 | 0.9917 | 0.9843 | 0.9618 | 0.9617 | 0.9583 | 0.9632 | 0.9716 | 0.9883 | 0.9929 | 0.9960 |
| 5/2 | 0.9759 | 0.9519 | 0.9220 | 0.9192 | 0.9237 | 0.9323 | 0.9389 | 0.9704 | 0.9835 | 0.9920 |
| 3 | 0.9449 | 0.9391 | 0.8631 | 0.8627 | 0.8876 | 0.8933 | 0.8981 | 0.9526 | 0.9708 | 0.9880 |
| 4 | 0.7917 | 0.8049 | 0.7443 | 0.7495 | 0.7736 | 0.7921 | 0.8164 | 0.8989 | 0.9392 | 0.9720 |
| 5 | 0.6323 | 0.6708 | 0.6047 | 0.6225 | 0.6582 | 0.6741 | 0.7103 | 0.8363 | 0.9012 | 0.9520 |

Table shows that the PDS with $d \leq 5/2$ are preferable than that of $d > 5/2$ for all range of $\lambda$. While this property of PDS more pronounced for the very sparse and dense models. It is surprise that for the moderate $\lambda$ the PDS with parameter $d \in [1/3, 2]$ appears to be asymptotically more correlated with chi-square statistic than the log-likelihood ratio statistic, where $d = 0$. But log-likelihood ratio statistic exhibit high limiting correlation with chi-square-statistic than the PDS with $d < 0$, i.e. satisfying Cramèr condition, $0.9335 \leq \rho(h_0, \lambda) \leq 1$ and $\arg\min \rho(h_0, \lambda) = 2.3750$. The PDS $CR_N(2/3)$ exhibit highest limiting correlation with chi-square-statistic for all range of $\lambda$: $0.9900 \leq \rho(h_{2/3}, \lambda) \leq 1$. This confirms recommendation of Cressie and Read (1984, p.462).

We close this section by indicating some of the remaining open problems in the study of intermediate properties of $h$-tests, the progress in solving of which depends on the progress in the theory of large deviations for respective test statistics.

- For very sparse models there remain a gap in the study of the $\alpha - IAE$ properties of $h_d$-tests for alternatives such that $\varepsilon(N) \geq (n^{d^*} \lambda_n^{d^*-1})^{-1/(2d^*+1)}$. Here we conjecture that instead of $d^* = \max(1, d)$, which actually is appeared due to condition (i) of Theorem 2.2 of Mirakhmedov (2020, see Remark 2.1), it should stay $d^* = \max(0, d)$. Then one would expect that equation (3.1) would be held for each $d \in (-1, 0]$ and family $\Im_{alt}$. In turn it would be possibly to extend the equality $e(h_{d_1}, h_{d_2}) = 1$ for any $d_1, d_2 \in (-1, 0]$ and family $\Im_{alt}$.

- For sparse and dense models when $\sqrt{n} \leq N \leq cn$ there remain a gap in the study of the properties of the chi-square test for alternatives such that $(n\lambda_n^2)^{-1/3} \leq \varepsilon(N)$
$= O\left((n\lambda_n)^{-1/3} \log^{2/3}(N^2/n)\right)$, and of the log-likelihood ratio test for alternatives such that $(n\lambda_n^2)^{-1/3} \leq \varepsilon(N) = O\left((n\lambda_n)^{-1/3}\right)$.

## 5. Proofs.

**Proof of Theorem 3.1, 3.2 and 3.3** follow from Assertions 1 and 2, 3 and 4, respectively, by putting $x_N = x_N(h) = \sqrt{n\lambda_n/2}\rho(h, \lambda_n)\varepsilon(N)$. The condition for the decreasing rate of $\varepsilon(N)$, which determines the corresponding family of alternatives, follows from the corresponding condition for the variable $x_N = x_N(h)$ of the assertions. For instance, when we apply Assertion 4 to prove Theorem 3.3 we have $x_N = \sqrt{n\lambda_n/2}\varepsilon(N) \ll N^{1/6}$ (since $\rho(h_d, \lambda_n) \to 1$ if $\lambda_n \to \infty$) which yield $\varepsilon(N) \ll (n\lambda_n^2)^{-1/3}$.

**Proof of Theorem 3.4**. For the chi-square statistic $h(u) = (u - \lambda_n)^2 / \lambda_n$, $Eh(\xi) = 1$, $\sigma^2(h) = 2$ and $\rho(\chi_N^2, \lambda_n) = 1$, hence

$$e_N^\alpha(\chi_N^2) = -\log P_0\left\{\chi_N^2 > x_N\sqrt{2N} + N\right\}, \text{ where } x_N = \sqrt{n\lambda_n/2}\,\varepsilon(N). \tag{5.1}$$

Under conditions of part (i) $x_N \to \infty$, $x_N = o\left((\sqrt{N}\min(1, \lambda_n^2))^{1/3}\right)$. Next, for each sequence of alternatives of $\Im_\gamma$ and $N \in \mathbb{N}_\gamma$ we have $N \ll \sqrt{n}$, $x_N \to \infty$, $x_N = o(\sqrt{N})$ and $N^{3/2}/n^{1/2}x_N \to 0$. Therefore, parts (i) and (ii) follows by applying in (5.1) Assertion 5.

*Proof of part* (iii). Set $v(n) = \lfloor \lambda_n + \sqrt{\lambda_n + n\lambda_n\varepsilon(N)} \rfloor + 1$. By (5.1) and definition of $\chi_N^2$ we have

$$P_0\left\{(\chi_N^2 - N)/\sqrt{2N} > x_N\right\} = P_0\left\{\sum_{m=1}^N \left((\eta_m - \lambda_n)^2 - \lambda_n\right) > n\lambda_n\varepsilon(N)\right\}$$

$$\geq P_0\left\{\sum_{m=2}^N \left((\eta_m - \lambda_n)^2 - \lambda_n\right) \geq 0 \,\Big/\, \eta_1 = v(n)\right\} P_0\left\{\eta_1 = v(n)\right\}$$

$$= P_0\left\{\sum_{m=1}^{N-1} \left((\bar{\eta}_m - \bar{\lambda}_n)^2 - \bar{\lambda}_n\right) \geq 0\right\} P_0\left\{\eta_1 = v(n)\right\}. \tag{5.2}$$

Here $\bar{\eta}_m \sim Bi\left(n - v(n), (N-1)^{-1}\right)$. Set $\bar{\lambda}_n = (n - v(n))/(N-1)$. It is easy to see that

$v(n)/n = \left(N^{-1} + \sqrt{\varepsilon(N)/N}\right)(1 + o(1))$ and $\bar{\lambda}_n = \lambda_n\left(1 + O\left(N^{-1} + \sqrt{\varepsilon(N)/N}\right)\right)$. We have

$$P_0\left\{\sum_{m=1}^{N-1}\left((\bar{\eta}_m - \lambda_n)^2 - \lambda_n\right) \geq 0\right\} \geq P_0\left\{\sum_{m=1}^{N-1}(\bar{\eta}_m - \bar{\lambda}_n)^2 \geq (N-1)\lambda_n\right\}$$

$$= P_0\left\{\sum_{m=1}^{N-1}\left((\bar{\eta}_m - \bar{\lambda}_n)^2 - \bar{\lambda}_n\right) \geq (v(n) - \lambda_n)\right\}$$

$$= P_0\left\{\frac{\sum_{m=1}^{N-1}\left((\bar{\eta}_m - \bar{\lambda}_n)^2 - \bar{\lambda}_n\right)}{\sqrt{2(n - v(n))^2/(N-1)}} \geq \sqrt{\varepsilon(N)}(1 + o(1))\right\} \geq c > 0, \qquad (5.3)$$

because $(n - v(n))\bar{\lambda}_n = n\lambda_n(1 + o(1)) \to \infty$, and hence the CLT for the statistic $\sum_{m=1}^{N-1}(\bar{\eta}_m - \bar{\lambda}_n)^2$ is enable to use, see Mirakhmedov (1992, Corollary 3).

Set $g(x, p) = x\log(x/p) + (1-x)\log((1-x)/(1-p))$, $x \in (0,1)$ and $p \in (0,1)$. Let $\varsigma \sim Bi(k, p)$. Due to Lemma 1 of Quine and Robinson (1985): for an integer $kx$

$$P\{\varsigma = kx\} \geq 0.8(2\pi kx(1-x))^{-1/2}\exp\{-kg(x, p)\}. \qquad (5.4)$$

Note that under $H_0$ the r.v. $\eta_1 \sim Bi(n, N^{-1})$, therefore applying (5.4) we obtain

$$P_0\{\eta_1 = v(n)\}$$

$$\geq c\left(v(n)(1 - v(n)n^{-1})\right)^{-1/2}\exp\left\{-v(n)\log(\lambda_n^{-1}v(n)) - n(1 - n^{-1}v(n))\log\frac{1 - v(n)n^{-1}}{1 - N^{-1}}\right\}$$

$$\geq c(v(n))^{-1/2}\exp\{-v(n)\log(\lambda_n^{-1}v(n))\}.$$

We have $n\varepsilon(N) \gg \sqrt{N}$, since (2.9), and $N\varepsilon(N) = (\omega_n N^2/n)^{2/3}$, $\log(v(n)/\lambda_n) < 2^{-1}\log(2N\varepsilon(N))$. Hence

$$-\frac{\log P_0\{\eta_1 = v(n)\}}{n\lambda_n\varepsilon^2(N)} \leq c\frac{\log v(n) + v(n)\log(\lambda_n^{-1}v(n))}{n\lambda_n\varepsilon^2(N)} \leq c\frac{\lambda_n + \sqrt{2\varepsilon(N)n\lambda_n}}{n\lambda_n\varepsilon^2(N)}\log\frac{v(n)}{\lambda_n}$$

$$\leq c\left(\frac{\log(N\varepsilon(N))}{n\varepsilon^2(N)} + \frac{\log(2\omega_n N^2/n)}{\sqrt{\varepsilon^3(N)n\lambda_n}}\right)$$

$$\leq \frac{c}{\sqrt{\varepsilon^3(N)n\lambda_n}}\left(\frac{\log(N\varepsilon(N))}{\sqrt{N\varepsilon(N)}} + \log\omega_n + \max(1, \log(N^2/n))\right) = o(1), \qquad (5.5)$$

since $\varepsilon^3(N)n\lambda_n = \omega_n^2$, for every sequence of intermediate alternatives of part (iii). Part (iii) follows from (5.1), (5.2), (5.3) and (5.5).

**Proof of Theorem 3.5**. Recall $\Lambda_N$ is PDS with parameter $d = 0$. Parts (i) and (ii) follows from Theorem 3.1 (i) and Theorems 3.2 (i) and 3.3, respectively. Further, for statistic $\Lambda_N$ we have $A_0(h) = 1$, $\sigma_0^2(h) = 2$ (see (2.4)) and $\rho(\Lambda_N, \lambda_n) \to 1$, if $\lambda_n \to \infty$, hence $x_N = \sqrt{n\lambda_n/2}\,\varepsilon(N)\rho(\Lambda_N,\lambda_n) = \sqrt{n\lambda_n/2}\,\varepsilon(N)(1+o(1))$ and $e_N^\alpha(\Lambda_N) = -\log P_0\{\Lambda_N > x_N\sqrt{2N}+N\}$.

Under the conditions of both (iii) and (iv) parts one can easily observe that $x_N = o(\sqrt{N})$ and $N^{3/2}/x_N^2\sqrt{n} \to 0$. Therefore, these parts follow from Assertion 6.

**Proof of Theorem 3.6** follows straightforwardly from Assertion 7.

**Acknowledgements**. The author thank anonymous reviewer for his valuable constructive comments that lead to an improved version of the paper.

**Appendix**. We still use notation of the previous sections. In addition we set: $\xi_m \sim Poi(np_m)$,

$$A_N(h) = N^{-1}\sum_{m=1}^N Eh(\xi_m), \qquad \tau_N(h) = \frac{1}{n}\sum_{m=1}^N \text{cov}(h(\xi_m),\xi_m),$$

$$\sigma_N^2(h) = N^{-1}\sum_{m=1}^N \text{var}\, h(\xi_m) - \lambda_n\tau_N^2(h), \quad \tilde{S}_N^h = (S_N^h - NA_N(h))/\sqrt{N\sigma_N^2(h)}. \tag{A.1}$$

**Assertion 1**. Let the function $h$ be not linear, $\lambda_n \to \lambda \in (0,\infty)$, $Np_{\max} \leq c_1$, some $c_1 > 0$,

$$\max_{1\leq m\leq N} E\exp\{H|h(\xi_m)|\} \leq c_2, \text{ and } \sigma_N^2(h) \geq c_3, \text{ some } H > 0, c_2 > 0 \text{ and } c_3 > 0.$$

Then for $x_N \to \infty$, $0 \leq x_N = o(N^{1/2})$ it holds

$$\log P\{\tilde{S}_N^h > x_N\} = -\frac{1}{2}x_N^2 + O\left(\log x_N + \frac{x_N^3}{\sqrt{N}}\right).$$

Assertion 1 follows from Theorem 2 of Ivchenko and Mirakhmedov (1995) and the fact that

$$1 - \Phi(x_n) = (x_n\sqrt{2\pi})^{-1}\exp\{-x_n^2/2\}(1+o(1)),\; x_n \to \infty. \tag{A.2}$$

**Assertion 2**. Let $\lambda_n \to \lambda \in (0,\infty)$ and $Np_m \leq c$. Then

$$\log P\{\tilde{S}_N^{h_d} > x_N\} = -\frac{1}{2}x_N^2(1+o(1)) \tag{A.3}$$

is valid if $x_N \to \infty$ and

(i) $x_N = o(\sqrt{N})$, for every $d \in (-1,0]$,

(ii) $x_N = o(N^{1/2(1+2d)})$, for every integer $d \geq 1$,

(iii) $x_N = o(\min(N^{1/8}, N^{1/2(1+2d^*)}))$, where $d^* = \max(1,d)$, for every non-integer $d > 0$.

**Assertion 3**. Let $\lambda_n \to 0$, $n\lambda_n^3 \to \infty$ and $p_1 = ... = p_N = N^{-1}$. Then

(i) For $x_N \to \infty$, $x_N = o\left((n\lambda_n^3)^{1/2(1+2d^*)}\right)$, where $d \in (-1,0)$ or integer $d \geq 1$,

(ii) For $x_N \to \infty$, $x_N = o\left(\min(n^{1/8},(n\lambda_n^3)^{1/2(1+2d^*)})\right)$, where non-integer $d > 0$, or $d = 0$,

. one has

$$\log P\{\tilde{S}_N^{h_d} > x_N\} = -\frac{1}{2}x_N^2(1+o(1))$$

$$\log P\{S_N^{h_d} > x_N |2^d - 1|\sqrt{2n\lambda_n} + n\} = -\frac{1}{2}x_N^2(1+o(1)).$$

(iii) For $x_N \to \infty$, $x_N = o\left(\min(n^{1/8},(n\lambda_n^3)^{1/6})\right)$ one has

$$\log P\{S_N^{h_0} > x_N\sqrt{2n\lambda_n \log 2} + n\lambda_n \log 2\} = -\frac{1}{2}x_N^2(1+o(1)).$$

**Assertion 4**. If $\lambda_n \to \infty$ and $Np_m \geq c > 0$, then for every $d > -1$ and $x_N \to \infty$, $x_N = o(N^{1/6})$

$$\log P\{S_N^{h_d} > x_N\sqrt{2N} + N\} = -\frac{1}{2}x_N^2(1+o(1)).$$

Assertions 2, 3 and 4 follow in immediate manner from Theorems 3.1, 3.2 and Theorems 3.5 and 3.8 of Mirakhmedov (2020), respectively, and (A.2).

**Assertion 5.** (i) Let $p_1 = ...p_N = N^{-1}$. Then for arbitrary $\lambda_n$ and $x_N \to \infty$, $x_N = o\left((\sqrt{N}\min(1,\lambda_n^2))^{1/3}\right)$ one has

$$\log P\{\chi_N^2 > x_N\sqrt{2N} + N\} = -\frac{1}{2}x_N^2(1+o(1)).$$

(ii) Let $N \ll \sqrt{n}$, $\min_m Np_m \geq c > 0$. If $x_N \to \infty$, $x_N = o(\sqrt{N})$ and $N^{-3/2}n^{1/2}x_N \to \infty$ then

$$\log P\{\chi_N^2 > x_N\sqrt{2N} + N\} = -\frac{1}{2}x_N^2 + O\left(\frac{x_N^3}{\sqrt{N}} + \log N + \frac{x_N N^{3/2}}{\sqrt{n}}\right);$$

**Assertion 6.** Let $\lambda_n \to \infty$ and $Np_{\min} \geq c > 0$. If $x_N \to \infty$ and $x_N = o(\sqrt{N})$ then

$$\log P\{\Lambda_N > x_N\sqrt{2N} + N\} = -\frac{1}{2}x_N^2 + O\left(\frac{x_N^3}{\sqrt{N}} + \log N + \frac{N^{3/2}}{\sqrt{n}}\right);$$

Assertions 5 (i) follows from Corollary 4.5 of Mirakhmedov (2020). Assertion 5 (ii) and Assertion 6 are Eq. (2.17) and Eq. (2.13), respectively, of Kallenberg (1985).

Set $\pi_r(\lambda) = \lambda^r e^{-\lambda}/r!$, $\sigma_r^2 = \pi_r(\lambda_n)(1-\pi_r(\lambda_n)) - \lambda_n^{-1}((r-\lambda_n)\pi_r(\lambda_n))^2$.

**Assertion 7**. Let $\lambda_n \to 0$, $n^{1/5}\lambda_n \to \infty$. Then for all $x_N \to \infty$, $x_N = o(N^{1/6}\lambda_n)$ and $r = 0,1,2$ one has

$$\log P\{\mu_r > x_N\sigma_r\sqrt{N} + N\pi_r(\lambda_n)\} = -\frac{1}{2}x_N^2(1+o(1)).$$

$$\log P\{C_n > x_N\sigma_0\sqrt{N} + N\pi_0(\lambda_n)\} = -\frac{1}{2}x_N^2(1+o(1))$$

Assertion 7 straightforwardly follows from Corollary 4.12 part (ii) of Mirakhmedov (2020) and fact that $C_n = \mu_0 - (n-N)$.

**References.**